\documentclass[a4paper,12pt]{article}
 \usepackage[top=2.5cm,bottom=2.5cm,left=2.5cm,right=2.5cm]{geometry}
 \usepackage{cite, amsmath, amssymb}
\newtheorem{prop}{Proposition}

 \title{\vspace*{2.5cm}{\bf Kemeny's constant and  Kirchhoffian indices for a family of non-regular graphs} }
\date{}
\author{{\bf Jos\'e Luis Palacios}\\
{\it Electrical and Computer Engineering Department,}\\{\it The University of New Mexico, Albuquerque, NM 87131, USA}\\
jpalacios@unm.edu\\
and\\
{\bf Greg Markowsky}\\
{\it Department of Mathematics,}\\
{\it Monash University, Melbourne, Australia}\\greg.markowsky@monash.edu.au}
\begin{document}
 \thispagestyle{empty} 
\maketitle


\begin{abstract} 
We find closed form formulas for Kemeny's constant and its relationship with two Kirchhoffian indices for some composite graphs that use as  basic building block a graph endowed with one of several symmetry properties.
\end{abstract}

\section{Introduction}

Let $G=(V,E)$ be a finite simple connected graph with vertex set $V=\{1, 2, \ldots, n\}$, edge set $E$ and  degrees $d_1\ge d_2\ge \cdots \ge d_n$. An automorphism of a graph is a bijection of $G$ onto $G$ that preserves adjacencies.  A graph is $d$-regular if all its vertices have degree $d$; it is vertex-transitive if there exists a graph automorphism that sends any vertex into any other vertex; it is edge-transitive if there exists a graph automorphism that sends any (undirected) edge into any other edge; finally it is distance regular if it is $d$-regular with diameter $D$ and there exist positive integers $b_0=d, b_1, \ldots, b_{D-1}$, $c_1=1, c_2, \ldots, c_D$ such that for every pair of vertices $u, v$ at distance $j$ apart, we have
\begin{itemize}
\item
the number of nodes at distance $j-1$ from $v$ which are neighbors of $u$ is $c_j$, $1\le j \le D$, and
\item
the number of nodes at distance $j+1$ from $v$ which are neighbors of $u$ is $b_j$, $0\le j \le D-1,$
\end{itemize}
and the numbers $c_j$ and $b_j$ are independent of the pair $u, v$ chosen.

A graph is walk-regular if for every $k$, the number of closed walks (starting and ending at $x$) of length $k$ is the same regardless of the vertex $x$, or equivalently, if every power $A^k$ of the incidence matrix $A$ of $G$ has the same values along the diagonal.
For these and all other concepts in Graph Theory not mentioned explicitly here we refer the reader to \cite{BH}.

The simple random walk on $G$ is the Markov chain $X_n, n\ge 1$ that jumps from one vertex of $G$ to a neighboring vertex with uniform probabilities.  If $P$ is the transition matrix of this chain, the stationary distribution of the random walk is the unique probabilistic vector $\pi$ that satisfies  $\pi P=\pi$ and that can be explicitly given as $\pi_i=\frac{d_i}{2|E|}$. The hitting time $T_b$ of the vertex $b$ is defined as $T_b=\inf \{n: X_n=b\}$ and its expectation, when the process starts at vertex $a$ is denoted by $E_aT_b$.  The Kemeny constant is defined as 
$$K=\sum_j \pi_j E_iT_j,$$
which turns out to be independent of $i$. For this fact and all other probabilistic notions we refer the reader to \cite{GS}.

The Kirchhoff index is a molecular descriptor defined on an undirected connected graph as
$$R(G)=\sum_{i<j} R_{ij},$$
where $R_{ij}$ is the effective resistance between $i$ and $j$ computed on the graph when it is thought of as an electric network with unit resistors on each edge. The multiplicative degree-Kirchhoff index is similarly defined as
$$R^*(G)=\sum_{i<j} d_i d_jR_{ij}.$$

The study of Kirchhoff indices and related quantities to do with electric resistances on graphs is an active research field. Recent works related to the results of this paper include  \cite{atik, ciardo, deville, Faught, geng, li, lei, zhou, zhou2}.

Chandra et al showed in \cite{C} that there is a close relationship between hitting times and effective resistances
\begin{equation}
\label{chandra}
E_aT_b+E_bT_a =2|E|R_{ab},
\end{equation}
and this can be used to express the Kirchhoff index in terms of hitting times (see \cite{P2001})
$$R(G)=\frac{1}{2|E|}\sum_{i<j}E_iT_j+E_jT_i.$$

In \cite{P1} and \cite{PR} it was noticed that for $d$-regular graphs there is a simple relationship between the Kirchhoff index and the Kemeny constant:
\begin{equation}
\label{simple}
R(G)=\frac{|V|}{d}K,
\end{equation}
but the relationship between $K$ and $R(G)$ is not that straightforward when the graph is not regular. There is a sustained interest in finding more closed form expressions or approximate computations for $K$, in theoretical and applied contexts.  We mention, for instance how in \cite{W} they found a general expression for $K$ in terms of the Moore-Penrose inverse of a relative of the Laplacian matrix, as well as upper and lower bounds for $K$ and some relations to the Kirchhoffian indices $R(G)$ and $R^*(G)$; also, in \cite{KD} they found some closed form expressions for $K$ in some cases of non-regular graphs as well as some approximations; finally, in \cite{X} they studied $K$ for some real-life scale-free networks and develop a randomized algorithm that approximately computes the Kemeny constant for any connected graph in nearly linear time with respect to the number of edges.

In this article we deal with a family of non-regular graphs built from certain symmetric graphs and find an explicit expression for $K$ of the composite graph in terms of $K$ of the building block graph, as well as the explicit relationship between $R(G)$ and $K$, which is more involved than (\ref{simple}).

\section{Highly symmetric graphs}
In \cite{DS} it was shown that for all vertex-transitive or distance regular graphs  the expected hitting times satisfy 
\begin{equation}
\label{main}
E_aT_b=E_bT_a~~~~~~{\rm for~all~pairs~}a, b {\rm ~in~} V.
\end{equation}

Graphs belonging to either of these families are called there ``highly symmetric (HS) graphs" and we will adopt here that terminology.   The set of HS graphs is not exhausted by those in \cite{DS}, i.e., vertex-transitive or distance-regular.   Indeed, these families of graphs are contained in the family of walk-regular graphs, which are HS, as was shown in \cite{G}.

In \cite{Pal98} we showed that if $G$ is edge-transitive then
(\ref{main}) holds
whenever the distance between $a$ and $b$ is even.  We also showed that
$$E_aT_b=E_bT_a=|V|-1,$$
if $d(a,b)=1$ and $G$ is edge-transitive and regular.  Based on these results we were interested in showing further  that an edge-transitive regular is HS, and perhaps that not all such graphs are walk-regular.  The second goal is futile, as there is a proof by Godsil in \cite{Go} showing that indeed all edge-transitive regular graphs are walk-regular.  However, we still believe that showing directly with electric arguments that regular edge-transitive graphs are HS is worthwhile, as an alternative to the algebraic graph theory route, providing different insights.
\vskip .2 in
For that purpose, we will be using a finer one-directional version of (\ref{chandra}),  shown in \cite{T},  stating that for all hitting times of simple random walk on $G$ we have
\begin{equation}
\label{tetali}
E_aT_b=\frac{1}{2}\sum_{z\in V} d_z \left(R_{ab}+R_{az}-R_{bz}\right).
\end{equation}
Our first result is a minor improvement of Theorem 2.1 in \cite{Pal98}
\begin{prop}
If there exists an automorphism $\Gamma$ such that $\Gamma(a)=b$ then $E_aT_b=E_bT_a$.
\end{prop}

{\bf Proof.} The automorphism ensures that $\sum_{z\in V}d_zR_{az}=\sum_{z\in V}d_zR_{bz}$.  Then clearly from (\ref{tetali})
$$E_aT_b=\frac{R_{ab}}{2}\sum_{z\in V}d_z=|E|R_{ab}=E_bT_a ~~~\bullet$$
The previous proposition clearly implies that a vertex-transitive graph is HS.  The next proposition is crucial for what follows.
\begin{prop}
Suppose $a,b,c,d\in V$, and there are automorphisms $\Phi, \Psi$ such that $\Phi(a)=c, \Psi(b)=d$. Then 
$$E_aT_b - E_bT_a = E_cT_d - E_dT_c.$$
\end{prop}

{\bf Proof.} By (\ref{tetali})  we have $E_aT_b - E_bT_a = \sum_{z\in V} d_zR_{zb} - \sum_{z \in V}d_z R_{za}$. The automorphisms ensure that $\sum_{z \in V}d_z R_{za} = \sum_{z \in V}d_z R_{zc}$ and $\sum_{z\in V} d_z R_{zb} = \sum_{z\in V} d_z R_{zd}$. Thus

$$
E_aT_b - E_bT_a = \sum_{z\in V} d_z R_{zb} - \sum_{z \in V}d_z R_{za} = \sum_{z \in V}d_z R_{zd} - \sum_{z \in V}d_z R_{zc} = E_cT_d - E_dT_c ~~~ \bullet
$$ 
\vskip .2 in
A remarkable characteristic of Proposition 2 is that we do not need to worry about where $\Phi$ takes $b$ or where $\Psi$ takes $a$; they essentially operate independently. We may make use of this proposition in the following way: define an equivalence relation on the vertices of $G$ by $a \equiv b$ whenever there is an automorphism taking $a$ to $b$. This equivalence relation partitions the vertices of $G$ into equivalence classes, let us say $A_1, \ldots, A_k$, and Proposition 2 shows that the quantities of the form $E_aT_b - E_bT_a$ depend only upon the equivalence classes of $a$ and $b$, not on the choices of $a$ and $b$ within their classes. In relation to the search for HS graphs, this means we need only to test whether hitting times are symmetric over a set of representatives, one from each equivalence class, rather than over all vertices in $G$. Note also that by Proposition 1, we do not need to test vertices within the same class. An immediate consequence is as follows.

\begin{prop}
Suppose $G$ is edge-transitive. Then $G$ is HS if and only if $G$ is regular.
\end{prop}

{\bf Proof:} Assume that $G$ has $q$ edges, one of which is $ab$.  There are at least $q$ automorphisms $\Gamma_1,  \ldots, \Gamma_q$ sending $ab$ onto the $q$ edges of $G$.  Define $V_1=\{\Gamma_1(a), \ldots, \Gamma_q(a)\}$ and $V_2=\{\Gamma_1(b), 
 \ldots, \Gamma_q(b)\}$.  Clearly $V_1\cup V_2=V$ and if $c, d\in V_1$ then $c=\Gamma_i(a)$ and $d=\Gamma_j(a)$ for some $i,j$, and thus $c=\Gamma_i(\Gamma_j^{-1}(d))$, implying by Proposition 1 that $E_cT_d=E_dT_c$.  A similar argument works if $c, d\in V_2$.  So we are only left to check the value of $E_cT_d-E_dT_c$ when $c\in V_1$, $d\in V_2$. By proposition 2 this value is equal to $E_aT_b-E_bT_a$ for the edge $ab$.
 
 But according to Theorem 2.4 in \cite{Pal98} we have $E_bT_a=\frac{2|E|}{d_1}-1$ and $E_aT_b=\frac{2|E|}{d_2}-1$ where $d_1$ (resp. $d_2$) is the common degree of all vertices in $V_1$ (resp. $V_2$) so that
 $$E_aT_b-E_bT_a=2|E|\left(\frac{1}{d_2}-\frac{1}{d_1}\right),$$
 which is equal to 0 if and only if the graph is regular~~$\bullet$
 
\vskip .2 in
It is well known (see \cite{BH}) that if $G$ is edge-transitive and regular, either $G$ is vertex-transitive (in which case all vertices belong to the same equivalence class, as discussed after Proposition 2) or bipartite, and the partition has two equivalence classes, $V_1$ and $V_2$ in the proof of proposition 3.  There are examples of regular edge-transitive graphs which are not vertex-transitive, the smallest of which is the Folkman graph,  depicted in \cite{F}, p. 227. If we insert a pair of axes X and Y, through the middle of the graph, one can see that the graph is not distance-regular because the top and bottom vertices on the X-axis, at distance 2, share exactly 2 neighbors, whereas the outermost pair of vertices on the Y-axis, also at distance 2, share 4 neighbors. This shows that edge-transitive regular graphs enlarge non-trivially the set of HS graphs.
\vskip .2 in

\section{Conjoining copies of an HS graph}

We will consider a building block graph $G_1=(V_1,E_1)$ that belongs to the HS class. For our purposes, the main properties with which $G_1$  is endowed are  that

(i) it is $d$-regular, and

(ii) $E_aT_b=E_bT_a=|E_1|R_{ab}$, on account of (\ref{chandra}).

\vskip .2 in
We can also compute the Kirchhoff index of any HS graph as follows: for any graph $G$ we have that
$$R(G)=\frac{1}{2}\sum_iR(i),$$
where
$$R(i)=\sum_j R_{ij}.$$
Now for a $d$-regular graph $G_1=(V_1,E_1)$ we have that
$$R(i)=\frac{1}{|E_1|}\sum_j E_iT_j=\frac{2}{nd}\sum_j E_iT_j=\frac{2}{d}\sum \pi_j E_iT_j =\frac{2}{d}K,$$
so we conclude that
 
 (iii) for an HS graph (in fact, it is enough that the graph be regular), all the quantities $R(i)$ have the same value, independent of $i$, and $R(G)=\frac{n}{2}R(1)$.
\vskip .2 in
Now let the graph $G=(E,V)$ consist of $\alpha$ copies of $G_1$, all of them conjoined, or glued,  at a single common vertex denoted by $c$.  This may be thought of as a model of a distributed network where each copy of $G_1$ is a ``province" and the node $c$ is the ``central   government". Then it is clear that all vertices of $G$ have degree $d$ except vertex $c$ that has degree $\alpha d$.  Also $|E|=\alpha |E_1|$ and $|V|=\alpha(|V_1|-1)+1$.  Let us denote by $K$ and $K_1$ the Kemeny's constants of $G$ and $G_1$, respectively.  With these preliminaries we can prove the following

\begin{prop}
For any HS graph $G_1$ and its composite graph $G$ defined above, and $\alpha\ge 1$,  we have
\begin{equation}
\label{A}
K=(2\alpha-1) K_1.
\end{equation}
\end{prop}

{\bf Proof.} We will compute $K$ as the sum of all expected hitting times starting from the vertex $c$ and normalized by the stationary distribution.  Since any vertex $i$ other than $c$ has degree $d$, it is clear that
$$\pi_i=\frac{d}{2|E|}=\frac{d}{2\alpha|E_1|}=\frac{1}{\alpha n}.$$
(Incidentally, $\pi_c=\frac{1}{n}$)

Now we have
$$K=\sum_{i\in V} \pi_i E_cT_i=\frac{1}{\alpha n}\sum_{i\in V}E_cT_i.$$

Now we will complicate a bit the notation and distinguish between $E_cT^G_i$ and $E_cT^{G_1}_i$, the hitting times starting from $c$ of the vertex $i$ when the random walk takes place on the whole $G$ and on the building block $G_1$, respectively.  Chandra et al.'s formula states that
$$E_iT^G_c+E_cT^G_i=2|E|R_{ic},$$
and there is no confusion where one computes the effective resistance $R_{ic}$, since it is has the same value whether it is computed in $G$ or in $G_1$. But clearly $E_iT^G_c=E_iT^{G_1}_c$, and we can solve 
$$E_cT^G_i=2|E|R_{ic}-E_iT^{G_1}_c=2|E|R_{ic}-|E_1|R_{ic}=(2\alpha-1)|E_1|R_{ic}.$$
Replacing above we see that
$$K=\frac{1}{\alpha n}\sum_{i\in V} (2\alpha-1)|E_1|R_{ic}=(2\alpha-1)\sum_{i\in V_1} \frac{1}{n}|E_1|R_{ic}=(2\alpha-1)\sum_{i\in V_1} \pi_iE_cT^{G_1}_i,$$
where $\pi_i$ is the stationary distribution of the walk on $G_1$, and we recognize the last summation as precisely $K_1$
$\bullet$

The example of the windmill graph $W(\eta, k)$ in \cite{KD} can be seen as the composite graph $G$ built from conjoining $\eta$ blocks equal to the complete graph $K_{k+1}$, whose Kemeny's constant is $K_1=\sum_i \pi_i E_cT_i=\frac{1}{k+1}\sum_i k=\frac{k^2}{k+1}$.  Then using our formula we see that
$$K=(2\eta-1)\frac{k^2}{k+1},$$
as stated in their Corollary 1.

Now we turn to the Kirchhoff index in the following

\vskip .2 in  
\begin{prop}
For any $n$-vertex HS $G_1$ and its composite graph $G$ defined above, and $\alpha\ge 1$,  we have
\begin{equation}
\label{B}
R(G)=\alpha\left[1+\frac{2(\alpha-1)}{n}\right]R(G_1).
\end{equation}
\end{prop}

{\bf Proof.} 

$$R(G)=\sum_{{i<j}\atop{i, j\in V}}R_{ij}=\alpha \sum_{{i<j}\atop{ i,j\in V_1} }R_{ij}+{\alpha \choose 2} \sum_{{i<j}\atop{ i\in V_1, j\in V_2}} R_{ij}.$$

Here the rightmost summation corresponds to all those effective resistances $R_{ij}$ between pairs of points on different copies $G_1$ and $G_2$ of the building block.  Since $c$ is a cut point, we have that $R_{ij}=R_{ic}+R_{cj}$ and thus
$$\sum_{{i<j}\atop {i\in V_1, j\in V_2}}R_{ij}=2\sum_i R_{ci}=\frac{2}{n}R(G_1),$$
where the last equality is due to property (iii). Putting all terms together we get the desired result $\bullet$

As a corollary, from the previous two propositions we get
\begin{prop} 
For the composite $G$ and $\alpha\ge 1$ we have
\begin{equation}
\label{C}
R(G)=\frac{\alpha(n+2\alpha-2)}{d(2\alpha-1)}K.
\end{equation}
\end{prop}

As another corollary, we can obtain the multiplicative degree-Kirchhoff index $R^*(G)$ in terms of the Kirchhoff index of the building block graph $G_1$ as follows:

\begin{prop} 
For the composite $G$ and $\alpha\ge 1$ we have
\begin{equation}
\label{D}
R^*(G)=\alpha(2\alpha-1) R^*(G_1)=\alpha(2\alpha-1)d^2R(G_1).
\end{equation}
\end{prop}

{\bf Proof.} For any graph $G$ we have that $R^*(G)=2|E|K$ (see \cite{PR}).  Also, if a graph is $d$-regular, (\ref{simple}) says that $R(G)=\frac{n}{d}K$.  Applying these two facts to the building block $G_1$ and its composite $G$ we get
$$R^*(G)=2|E|K=2\alpha|E_1|(2\alpha-1)K_1=\alpha(2\alpha-1)ndK_1=\alpha(2\alpha-1)d^2\frac{n}{d}K_1$$
$$=\alpha(2\alpha-1)d^2R(G_1)=\alpha(2\alpha-1)R^*(G_1) \bullet$$

Using (\ref{C}) and the fact that $R^{*}(G)=2|E|K$,  we get this relationship between the two Kirchhoffian indices of the composite graph:

\begin{prop} 
For the composite $G$ and $\alpha\ge 1$ we have
\begin{equation}
\label{E}
R^*(G)=\frac{(2\alpha-1)nd^2}{n+2\alpha-2}R(G).
\end{equation}
\end{prop}

\end{document}